# DUAL HYPERBOLIC CONCHOIDAL MOTION

**Mustafa Kazaz[1], Hasan Hüseyin UĞURLU[2]**

[1] Celal Bayar University, Faculty of Art and Science, Department of Mathematics, Muradiye Campus, Manisa, Turkey. E-mail: **mustafa.kazaz@cbu.edu.tr**

[2] Gazi University, Faculty of Education, Department of Secondary Education, Science and Mathematics Teaching, Mathematics Teaching Program, Ankara, Turkey. E-mail: **hugurlu@gazi.edu.tr**

**Abstract**

Dual spherical conchoidal motion has been defined by Yapar [6]. In this work, we define this motion on a dual hyperbolic unit sphere $\tilde{\mathbb{H}}_0^2$ in the dual Lorentzian space $\mathbb{D}_1^3$ with dual signature $(+, +, -)$, and the results carried to the Lorentzian lines space $IR_1^3$ by means of the Study's mapping. We also obtain the study maps of the orbits drawn on the fixed dual hyperbolic unit sphere by unit dual vectors of an orthonormal base $\{\tilde{v}_1, \tilde{v}_2, \tilde{v}_3\}$.

**2010 MR Subject Classification:** 53A17, 53B30, 53A40.

## 1. Introduction

Clifford (1873) introduced dual numbers in the form $\lambda + \varepsilon \lambda^*$ with $\varepsilon^2 = 0$ as a tool for his geometrical investigations, in particularly, for studying the Non-Euclidean geometry. After him, E. Study (1903) defined dual numbers as dual angles to specify the relations between two lines in the Euclidean space $\mathbb{E}^3$, and then he used dual numbers and dual vectors in his research on the geometry of lines and kinematics, and defined the mapping which is called with his name (E. Study's Mapping): There is a one to one correspondence between an oriented straight line in the Euclidean 3-space $\mathbb{E}^3$ and a dual point on the surface of a dual unit sphere $\mathbb{S}^2$ in the dual space $\mathbb{D}^3$ [3]. Then, a differentiable curve on the sphere $\mathbb{S}^2$ corresponds to a ruled surface in the line space $IR^3$ ([1], [5]).

Ruled surfaces have been widely applied in surface design and simulation of rigid bodies [8].

It is known that dual vectors, dual inner product, dual cross product, dual angle, dual orthogonal matrices, E. Study Mapping, etc. are the most important notions for applications of dual geometry to engineering. For example, the dual angle $\tilde{\varphi} = \varphi + \varepsilon \varphi^*$ between two dual unit vectors is formed with real angle $\varphi$ between corresponding two directed lines in the line space $IR^3$ and the shortest distance $\varphi^*$ between these directed lines. These notions lay the foundations for the study of spherical and spatial motions. Dual Lorentzian correspondences of these notions were introduced and also, several important theorems and results related to geometry of this space were given by the authors [4].

E. Study mapping plays a fundamental role between the real and dual Lorentzian spaces [8]. By this mapping, a curve on a dual hyperbolic unit sphere $\tilde{\mathbb{H}}_0^2$ corresponds to a timelike ruled surface in the Lorentzian line space $IR_1^3$, that is, there exist a one-to-one correspondence between the geometry of curves on $\tilde{\mathbb{H}}_0^2$ and the geometry of timelike ruled surfaces in $IR_1^3$. Similarly, a timelike (spacelike) curve on a dual Lorentzian unit sphere $\tilde{S}_1^2$ corresponds to a spacelike (or timelike) ruled surface in the Lorentzian line space $IR_1^3$, that is, there exist one-to-one correspondence between the geometry of timelike (spacelike) curves on $\tilde{S}_1^2$ and the geometry of spacelike (timelike) ruled surfaces in $IR_1^3$ [7]. Since the dual Lorentzian metric is indefinite, the angle concept in this space is very interesting. For instance, the dual hyperbolic angle $\tilde{\varphi} = \varphi + \varepsilon \varphi^*$ between two dual timelike unit vectors is a dual value formed with the (real) hyperbolic angle $\varphi$ between corresponding two directed timelike lines in the Lorentzian line space $IR_1^3$ and the shortest Lorentzian distance $\varphi^*$ between these directed timelike lines (Fig. 1.1. (a), (b)).

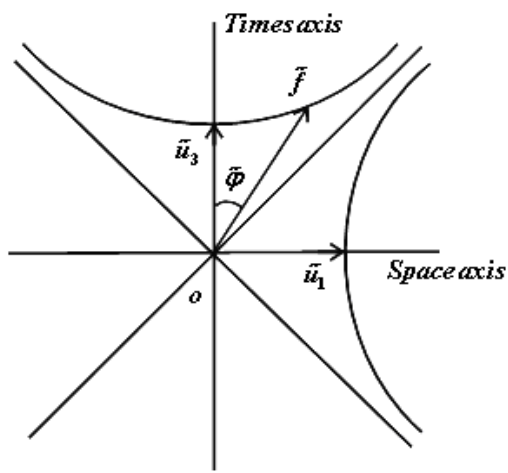 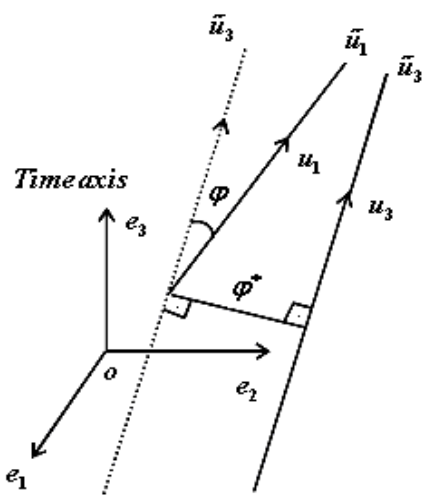

Fig.1.1. (a) The hyperbolic angle $\tilde{\varphi} = \varphi + \varepsilon \varphi^*$ between dual timelike unit vectors $\tilde{u}_3$ and $\tilde{f}$.

Fig.1.1. (b) The Lorentzian geometrical interpretation of the dual hyperbolic angle $\tilde{\varphi}$.

Planar and spherical conchoidal motions have been introduced by Karger and Novac [2]. Also, Yapar [6] defined the dual spherical conchoidal motion and gave the geometrical interpretations in the Euclidean lines space $IR^3$. In this paper, we define the dual hyperbolic analogy of planar, spherical and dual spherical conchoidal motions.

The Lorentzian motions in the Minkowski 3-space $IR_1^3$ are represented in the dual Lorentzian 3-space $\mathbb{D}_1^3$ by dual Lorentzian orthogonal $3 \times 3$ matrices $A = (a_{ij})$, $A^{-1} = \varepsilon A^T \varepsilon$, where $a_{ij}$ are dual functions of one dual variable $t$, and

$$\varepsilon = \begin{pmatrix} 1 & 0 & 0 \\ 0 & 1 & 0 \\ 0 & 0 & -1 \end{pmatrix}$$

is a signature matrix [4]. This means that when a Lorentzian motion is given in $IR_1^3$, we can find a corresponding dual Lorentzian orthogonal $3 \times 3$ matrix $A$.

## 2. Dual Lorentzian Vectors and E. Study Mapping

In this section, we give a brief summary of the theory of dual numbers and dual Lorentzian vectors. Let $IR_1^3$ be a 3-dimensional Minkowski space over the field of real numbers $IR$ with the Lorentzian inner product $<,>$ given by

$$<a,b> = a_1 b_1 + a_2 b_2 - a_3 b_3,$$

where $a = (a_1, a_2, a_3)$ and $b = (b_1, b_2, b_3) \in IR^3$. A vector $a = (a_1, a_2, a_3)$ of $IR_1^3$ is said to be *timelike* if $<a,a> < 0$, *spacelike* if $<a,a> > 0$ or $a = 0$, and *lightlike* (*null*) if $<a,a> = 0$ and $a \neq 0$. The norm of a vector $a$ is defined by $|a| = \sqrt{|<a,a>|}$. Now, let $a = (a_1, a_2, a_3)$ and $b = (b_1, b_2, b_3)$ be two vectors in $IR_1^3$, then the *Lorentzian cross product* of $a$ and $b$ is given by

$$a \times b = (a_3 b_2 - a_2 b_3, \, a_1 b_3 - a_3 b_1, \, a_1 b_2 - a_2 b_1).$$

A *dual number* has the form $\lambda + \varepsilon \lambda^*$, where $\lambda$ and $\lambda^*$ are real numbers and $\varepsilon$ stands for the dual unit which is subject to the rules:

$$\varepsilon \neq 0, \, \varepsilon^2 = 0, \, 0\varepsilon = \varepsilon 0 = 0, \, 1\varepsilon = \varepsilon 1 = \varepsilon.$$

Like a real number which can be considered as an angle, in differential geometry and motion analysis of spatial mechanisms, a dual number is also commonly referred as a dual angle $\lambda + \varepsilon \lambda^*$ between two lines in the space. The real part $\lambda$ of the dual angle is the projected angle between the lines, and the dual part $\lambda^*$ is the length along the common normal of the lines.

We denote the set of all dual numbers by $\mathbb{D}$:

$$\mathbb{D} = \{ \lambda + \varepsilon \lambda^* : \lambda, \lambda^* \in IR, \varepsilon^2 = 0 \}.$$

Equality, addition and multiplication are defined in $\mathbb{D}$ by

$$\lambda + \varepsilon \lambda^* = \beta + \varepsilon \beta^* \text{ iff } \lambda = \beta \text{ and } \lambda^* = \beta^*,$$

$$(\lambda + \varepsilon \lambda^*) + (\beta + \varepsilon \beta^*) = (\lambda + \beta) + \varepsilon (\lambda^* + \beta^*),$$

and

$$(\lambda + \varepsilon \lambda^*)(\beta + \varepsilon \beta^*) = \lambda \beta + \varepsilon (\lambda \beta^* + \lambda^* \beta),$$

respectively. Then it is easy to show that $(\mathbb{D}, +, .)$ is a commutative ring with unity. Moreover, if $\lambda + \varepsilon \lambda^*$, $\beta + \varepsilon \beta^* \in \mathbb{D}$ with $\beta \neq 0$ then the division is given by

$$\frac{\lambda + \varepsilon \lambda^*}{\beta + \varepsilon \beta^*} = \frac{\lambda}{\beta} + \varepsilon \left( \frac{\lambda^*}{\beta} - \frac{\lambda \beta^*}{\beta^2} \right).$$

Now, let $f$ be a differentiable function with dual variable $x + \varepsilon x^*$. Then the Maclaurin series generated by $f$ is

$$f(x + \varepsilon x^*) = f(x) + \varepsilon x^* f'(x),$$

where $f'(x)$ is the derivative of $f$.

Let $\mathbb{D}^3$ be the set of all triples of dual numbers, i.e.

$$\mathbb{D}^3 = \{ \tilde{a} = (a_1, a_2, a_3) \mid a_i \in \mathbb{D}, 1 \leq i \leq 3 \}.$$

The elements of $\mathbb{D}^3$ are called as *dual vectors*. A dual vector $\tilde{a}$ may be expressed in the form $\tilde{a} = a + \varepsilon a^*$, where $a$ and $a^*$ are the vectors of $IR^3$.

Now let $\tilde{a} = a + \varepsilon a^*$, $\tilde{b} = b + \varepsilon b^* \in \mathbb{D}^3$ and $\lambda = \lambda_1 + \varepsilon \lambda_1^* \in \mathbb{D}$. Then we define

$$\tilde{a} + \tilde{b} = a + b + \varepsilon (a^* + b^*),$$
$$\lambda \tilde{a} = \lambda_1 a + \varepsilon (\lambda_1 a^* + \lambda_1^* a).$$

Then $\mathbb{D}^3$ becomes a unitary $\mathbb{D}$-module with these operations. It is called $\mathbb{D}$-*module* or *dual space*.

The *Lorentzian inner product* of two dual vectors $\tilde{a}=a+\varepsilon a^*$, $\tilde{b}= b + \varepsilon b^*$ is defined by

$$<\tilde{a},\tilde{b}>=<a,b>+\varepsilon(<a,b^*>+<a^*,b>),$$

where $<a,b>$ is the Lorentzian inner product of the vectors $a$ and $b$ in the Minkowski 3-space $IR_1^3$. Then a dual vector $\tilde{a}=a+\varepsilon a^*$ is said to be *timelike* if $a$ is timelike, *spacelike* if $a$ is spacelike or $a=0$ and *lightlike* (*null*) if $a$ is lightlike (null) and $a \neq 0$. The set of all dual Lorentzian vectors is called *dual Lorentzian space* and it is denoted by $\mathbb{D}_1^3$:

$$\mathbb{D}_1^3 = \left\{ \tilde{a}=a+\varepsilon a^*: a, a^* \in IR_1^3 \right\}.$$

The *Lorentzian cross product* of dual vectors $\tilde{a}$ and $\tilde{b} \in \mathbb{D}_1^3$ is defined by

$$\tilde{a}\times\tilde{b} = a\times b+\varepsilon(a^*\times b+a\times b^*),$$

where $a \times b$ is the Lorentzian cross product in $IR_1^3$.

**Lemma 2.1.** *Let* $\tilde{a}, \tilde{b}, \tilde{c}, \tilde{d} \in \mathbb{D}_1^3$. *Then we have*

$$\tilde{a}\times\tilde{b}=-\tilde{b}\times\tilde{a},$$
$$<\tilde{a}\times\tilde{b}, \tilde{a}>=0 \; ; and \; <\tilde{a}\times\tilde{b}, \tilde{b}>=0,$$
$$<\tilde{a}\times\tilde{b},\tilde{c}>=-\det(\tilde{a}, \tilde{b}, \tilde{c})$$
$$(\tilde{a}\times\tilde{b})\times\tilde{c}=-<\tilde{a},\tilde{c}>\tilde{b}+<\tilde{b},\tilde{c}>\tilde{a},$$
$$<\tilde{a}\times\tilde{b}, \tilde{c}\times\tilde{d}>=-<\tilde{a},\tilde{c}><\tilde{b},\tilde{d}>+<\tilde{a},\tilde{d}><\tilde{b},\tilde{c}>. \qquad [4]$$

Let $\tilde{a}=a+\varepsilon a^* \in \mathbb{D}_1^3$. Then $\tilde{a}$ is said to be *unit dual timelike vector* (resp., *unit dual spacelike vector*) if the vectors $a$ and $a^*$ satisfy the following equations:

$$<a,a>=-1 (\text{resp.}, <a,a>=1), \; <a,a^*>=0.$$

The set of all unit dual timelike vectors (resp., all unit dual spacelike vectors) is called the *dual hyperbolic unit sphere* (resp. *dual Lorentzian unit sphere*), and is denoted by $\tilde{\mathbb{H}}_0^2$ (resp., $\tilde{\mathbb{S}}_1^2$) [4].

**Theorem 2.2** *(E. Study's Mapping) The dual timelike unit vectors of the dual hyperbolic unit sphere $\tilde{\mathbb{H}}_0^2$ are in one-to-one correspondence with the directed timelike lines of the Minkowski 3-space $IR_1^3$* [4].

## 3. Conchoidal Motion on the Dual Hyperbolic Unit Sphere $\tilde{H}_0^2$

We will define conchoidal motion on the dual hyperbolic unit sphere. Let us consider a fixed dual orthonormal frame $R = \{0; \tilde{u}_1, \tilde{u}_2, \tilde{u}_3 \text{(timelike)}\}$ and denote this frame by the dual hyperbolic unit sphere $H'$. Let $H_0^1$ be a great hyperbolic circle on $H'$, and $C$ be a point not lying on $H_0^1$. The frame $\{0; \tilde{u}_1, \tilde{u}_2, \tilde{u}_3\}$ is chosen as shown in Fig. 3.1, where $\tilde{u}_2$ and $\tilde{u}_3$ lie in the timelike plane of circle $H_0^1$, and the timelike plane $\tilde{u}_1\tilde{u}_3$ contains the chosen point $C$.

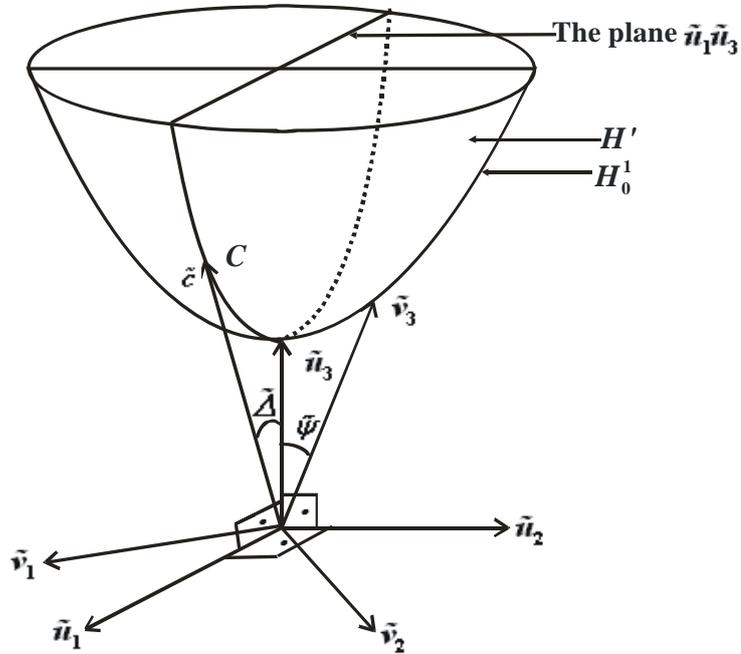

**Figure 3.1.** Dual hyperbolic conchoidal motion ($\tilde{c}$, $\tilde{u}_3$ and $\tilde{v}_3$ are dual timelike unit vectors and the others are unit vectors of dual Lorentzian unit sphere $\tilde{S}_1^2$).

Let us consider an orthonormal dual frame $\{0; \tilde{v}_1, \tilde{v}_2, \tilde{v}_3 \text{(timelike)}\}$ as shown in Fig.3.1. The frame $\{0; \tilde{v}_1, \tilde{v}_2, \tilde{v}_3\}$ moves now in such a way that the timelike vector $\tilde{v}_3$ rotates in great hyperbolic circle $H_0^1$ while the timelike plane $\tilde{v}_1\tilde{v}_3$ passes through the point $C$ all the time. As the parameter of the motion we choose the dual hyperbolic angle $\tilde{\psi} = \psi + \varepsilon\psi^*$ of the timelike vectors $\tilde{v}_3$ and $\tilde{u}_3$, where $\tilde{v}_3 = \tilde{u}_3 \cosh\tilde{\psi} + \tilde{u}_2 \sinh\tilde{\psi}$.

Further, let us denote

$$\tilde{v}_1 = A_1\tilde{u}_1 + A_2\tilde{u}_2 + A_3\tilde{u}_3, \quad \|\tilde{v}_1\| = A_1^2 + A_2^2 - A_3^2 = 1, \tag{1}$$

where $A_i$'s are dual numbers. By orthonormality, we have $<\tilde{v}_3, \tilde{v}_1> = 0$, i.e.

$$<\tilde{u}_3\cosh\tilde{\psi} + \tilde{u}_2\sinh\tilde{\psi},\ A_1\tilde{u}_1 + A_2\tilde{u}_2 + A_3\tilde{u}_3> = 0,$$

or

$$-A_3\cosh\tilde{\psi} + A_2\sinh\tilde{\psi} = 0. \tag{2}$$

Further, we may write $\tilde{c}$ as follows:

$$\tilde{c} = \tilde{u}_3\cosh\tilde{\Delta} + \tilde{u}_1\sinh\tilde{\Delta},$$

where $\tilde{\Delta} = \Delta + \varepsilon\Delta^*$ is dual hyperbolic angle between the timelike vectors $\tilde{c}$ and $\tilde{u}_3$. Since the timelike plane $\tilde{v}_1\tilde{v}_3$ has to pass through the point $C$ all the time, the vectors $\tilde{v}_1$, $\tilde{v}_3$ and $\tilde{c}$ must be co-planar, that is $\det(\tilde{v}_1, \tilde{v}_3, \tilde{c}) = 0$. Thus we get the equation

$$A_1\sinh\tilde{\psi}\cosh\tilde{\Delta} + A_2\cosh\tilde{\psi}\sinh\tilde{\Delta} - A_3\sinh\tilde{\psi}\sinh\tilde{\Delta} = 0. \tag{3}$$

Then we have three equations altogether for the unknown $A_1$, $A_2$, $A_3$:

$$A_1^2 + A_2^2 - A_3^2 = 1,$$
$$-A_3\cosh\tilde{\psi} + A_2\sinh\tilde{\psi} = 0,$$
$$(A_2\cosh\tilde{\psi} - A_3\sinh\tilde{\psi})\sinh\tilde{\Delta} + A_1\sinh\tilde{\psi}\cosh\tilde{\Delta} = 0.$$

From the second equation we obtain $A_3 = \tilde{\lambda}\sinh\tilde{\psi}$, $A_2 = \tilde{\lambda}\cosh\tilde{\psi}$, $\tilde{\lambda} \in \mathbb{D}$. Substituting this into the third equation, we obtain

$$\tilde{\lambda}\sinh\tilde{\Delta} + A_1\sinh\tilde{\psi}\cosh\tilde{\Delta} = 0,\ \text{i.e.,}\ A_1\sinh\tilde{\psi} = -\tilde{\lambda}\tanh\tilde{\Delta}.$$

Multiply the first equation by $\sinh^2\tilde{\psi}$, upon substitution we have

$$\sinh^2\tilde{\psi} = \tilde{\lambda}^2\left(\tanh^2\tilde{\Delta} + \sinh^2\tilde{\psi}\right),$$

hence

$$\tilde{\lambda} = \pm\sinh\tilde{\psi}\left(\sinh^2\tilde{\psi} + \tanh^2\tilde{\Delta}\right)^{-1/2}.$$

Choose the plus sign. Then consequently we have

$$\tilde{v}_1 = \left[-\left(\sinh^2\tilde{\psi} + \tanh^2\tilde{\Delta}\right)^{-1/2}\tanh\tilde{\Delta},\ \sinh\tilde{\psi}\cosh\tilde{\psi}\left(\sinh^2\tilde{\psi} + \tanh^2\tilde{\Delta}\right)^{-1/2},\right.$$
$$\left.\sinh^2\tilde{\psi}\left(\sinh^2\tilde{\psi} + \tanh^2\tilde{\Delta}\right)^{-1/2}\right],$$

$$\tilde{v}_3 = [0,\ \sinh\tilde{\psi},\ \cosh\tilde{\psi}],$$

and from $\tilde{v}_2 = \tilde{v}_1 \times \tilde{v}_3$ we obtain

$$\tilde{v}_2 = \left[ -\sinh\tilde{\psi}\left(\sinh^2\tilde{\psi} + \tanh^2\tilde{\Delta}\right)^{-\frac{1}{2}},\ -\cosh\tilde{\psi}\ \tanh\tilde{\Delta}\left(\sinh^2\tilde{\psi} + \tanh^2\tilde{\Delta}\right)^{-\frac{1}{2}},\right.$$
$$\left. -\sinh\tilde{\psi}\ \tanh\tilde{\Delta}\left(\sinh^2\tilde{\psi} + \tanh^2\tilde{\Delta}\right)^{-\frac{1}{2}} \right].$$

Thus a moving orthonormal dual frame $\{0;\ \tilde{v}_1,\ \tilde{v}_2,\ \tilde{v}_3\}$ is chosen. Let us represent this moving frame by dual hyperbolic unit sphere $H$. Then, a dual hyperbolic conchoidal motion which is analogous to the real conchoidal motion [2] is obtained. In this case, dual hyperbolic conchoidal motion is represented by $H/H'$.

Now, let us choose a fixed point $X$ on the trace of $H$ in the plane $\tilde{v}_1\tilde{v}_3$ (we should note that the trace of a surface in any plane is simply the intersection of the surface and the plane). During the dual hyperbolic conchoidal motion, the dual point $X$ draws an orbit on $H'$. We denote the dual hyperbolic angles of $\tilde{v}_1\tilde{x}$ and $\tilde{x}\tilde{v}_3$ by $P = p + \varepsilon p^*$ and $Q = q + \varepsilon q^*$, respectively (see Fig. 3.2).

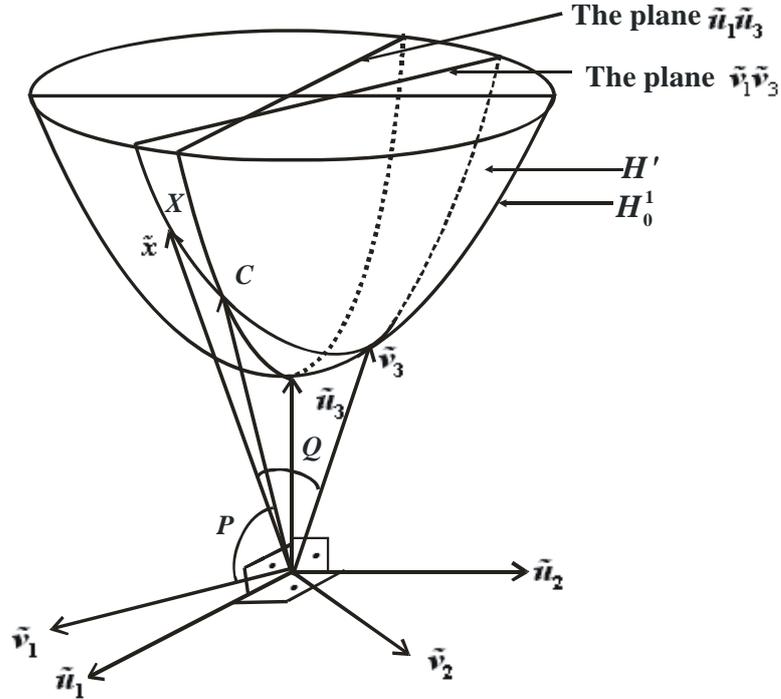

Figure 3.2: The timelike vector $\tilde{x}$ is on the plane $\tilde{v}_1\tilde{v}_3$.

Then we may write

$$\tilde{x} = \frac{\tilde{v}_1 \sinh P + \tilde{v}_3 \sinh Q}{\sinh(P+Q)}.$$

Since $P+Q = \frac{\pi}{2}$ and $\sinh(\pi/2) \neq 0$ we can write $\frac{1}{\sinh(P+Q)} = a = \text{constant}$. So, $\tilde{x}$ can be written as follows,

$$\tilde{x} = a(\tilde{v}_1 \sinh P + \tilde{v}_3 \sinh Q). \tag{4}$$

Making the necessary calculations for $\tilde{x}$ we have

$$x = a\left(A^{-1/2} \tanh \sigma \sinh p, \ -A^{-1/2} \sinh \psi \cosh \psi \sinh p + \sinh \psi \sinh q, \right.$$
$$\left. A^{-1/2} \sinh^2 \psi \sinh p + \cosh \psi \sinh q\right); \tag{5}$$

$$x^* = a\left(p^* A^{-1/2} \tanh \sigma \cosh p - \psi^* A^{-3/2} \tanh \sigma \sinh \psi \cosh \psi \sinh p - \sigma^* A^{-3/2} \tanh^2 \sigma \right.$$
$$\text{sech}^2 \sigma \sinh p + \sigma^* A^{-1/2} \text{sech}^2 \sigma \sinh p, \ -p^* A^{-1/2} \sinh \psi \cosh \psi \cosh p$$
$$+ q^* \sinh \psi \cosh q + \psi^* \cosh \psi \sinh q + \psi^* A^{-3/2} \sinh^2 \psi \cosh^2 \psi \sinh p$$
$$- \psi^* A^{-1/2} \cosh 2\psi \sinh p + \sigma^* A^{-3/2} \tanh \sigma \text{sech}^2 \sigma \sinh \psi \cosh \psi \sinh p,$$
$$\psi^* \sinh \psi \sinh q + q^* \cosh \psi \cosh q + A^{-3/2} \sigma^* \tanh \sigma \text{sech}^2 \sigma \sinh^2 \psi \sinh p$$
$$\left. - \psi^* A^{-3/2} \sinh^3 \psi \cosh \psi \sinh p + \psi^* A^{-1/2} \sinh 2\psi \sinh p + A^{-1/2} p^* \sinh^2 \psi \cosh p\right), \tag{6}$$

where $x$ and $x^*$ are the real and dual parts of $\tilde{x}$, respectively, $A = \sinh^2 \psi + \tanh^2 \psi$, $\Delta = \sigma + \varepsilon \sigma^* = \text{const}$. $P, Q = \text{const}$. Equations (5) and (6) depend on only two parameters $\psi$ and $\psi^*$. Thus, equations (5) and (6) represent a timelike line congruence in $IR_1^3$.

A timelike line congruence may be expressed as follows:

Let $m = m(\psi, \psi^*)$ be a position vector of the reference surface of a timelike line congruence, and let $x = x(\psi, \psi^*)$ be a unit vector in direction of a timelike line $x$ of the timelike line congruence. Let $y$ denote the position vector of an arbitrary point $Y = (y_1, y_2, y_3)$ of the fixed timelike line $x$ of the timelike line congruence in $IR_1^3$. Then we have $y = m + px$. We know that the moment vector $x^*$ of the vector $x$ with respect to the origin 0 is $x^* = m \times x$ and $x \times x^* = m + <m, x> x$. Then we may write $\lambda = p - \langle m, x \rangle$. Thus, we have

$$y = x(\psi, \psi^*) \times x^*(\psi, \psi^*) + \lambda x(\psi, \psi^*). \tag{7}$$

Since $(y_1, y_2, y_3)$ are the coordinates of $y$ we have

$$y_1 = a^2 \left( q^* A^{-1/2} \sinh\psi \cosh 2\psi \sinh p \cosh q - \psi^* A^{-1/2} \cosh\psi \cosh 2\psi \sinh p \sinh q \right.$$
$$+ \psi^* A^{-1} \sinh^2 \psi \sinh^2 p + \psi^* \sinh^2 q + \psi^* A^{-3/2} \sinh^2 \psi \cosh 2\psi \sinh p \sinh q \cosh p$$
$$+ \sigma^* A^{-3/2} \sinh\psi \tanh\sigma \operatorname{sech}^2 \sigma \sinh p \sinh q + 2\sigma^* A^{-2} \sinh^3 \psi \cosh\psi \tanh\sigma \quad (8)$$
$$\operatorname{sech}^2 \sigma \sinh^2 p - p^* A^{-1/2} \sinh\psi \cosh 2\psi \cosh p \sinh q \Big) + a\lambda A^{-1/2} \tanh\sigma \sinh p ;$$

$$y_2 = a^2 \left( q^* A^{-1/2} \tanh\sigma \cosh\psi \sinh p \cosh q - p^* A^{-1/2} \tanh\sigma \cosh\psi \sinh q \cosh p \right.$$
$$- \sigma^* A^{-1} \operatorname{sech}^2 \sigma \sinh^2 \psi \sinh^2 p + \sigma^* A^{-3/2} \tanh^2 \sigma \operatorname{sech}^2 \sigma \cosh\psi \sinh p \sinh q$$
$$- \sigma^* A^{-1/2} \operatorname{sech}^2 \sigma \cosh\psi \sinh p \sinh q + \psi^* A^{-1/2} \tanh\sigma \sinh\psi \sinh p \sinh q \quad (9)$$
$$+ \psi^* A^{-3/2} \tanh\sigma \sinh\psi \cosh^2 \psi \sinh p \sinh q + \psi^* A^{-1} \tanh\sigma \sinh 2\psi \sinh^2 p \Big)$$
$$+ a\lambda(-A^{-1/2} \sinh\psi \cosh\psi \sinh p + \sinh\psi \sinh q);$$

and

$$y_3 = a^2 \left( q^* A^{-1/2} \tanh\sigma \sinh\psi \sinh p \cosh q + \psi^* A^{-1/2} \tanh\sigma \cosh\psi \sinh q \sinh p \right.$$
$$- \psi^* A^{-1} \tanh\sigma \cosh 2\psi \sinh^2 p + \sigma^* A^{-1} \sinh\psi \cosh\psi \operatorname{sech}^2 \sigma \sinh^2 p$$
$$- p^* A^{-1/2} \tanh\sigma \sinh\psi \sinh q \cosh p - \sigma^* A^{-1/2} \sinh\psi \sinh p \sinh q \operatorname{sech}^2 \sigma \quad (10)$$
$$+ \psi^* A^{-3/2} \tanh\sigma \sinh^2 \psi \cosh\psi \sinh p \sinh q + \sigma^* A^{-3/2} \sinh\psi \tanh^2 \sigma$$
$$\operatorname{sech}^2 \sigma \sinh p \sinh q \Big) + a\lambda (A^{-1/2} \sinh^2 \psi \sinh p + \cosh\psi \sinh q).$$

If we take $Q = q + \varepsilon q^* = 0$, then the condition (4) gives us $X = \tilde{v}_1$. Thus, from equations (8), (9) and (10) we have

$$y_1 = \left( \psi^* \sinh^2 \psi (\sinh^2 \psi + \tanh^2 \sigma)^{-1} + 2\sigma^* \sinh^3 \psi \cosh\psi (\sinh^2 \psi + \tanh^2 \sigma)^{-2} \right)$$
$$+ \lambda \tanh\sigma (\sinh^2 \psi + \tanh^2 \sigma)^{-1/2}, \quad (11)$$

$$y_2 = \left( -\sigma^* \operatorname{sech}^2 \sigma \sinh^2 \psi + \psi^* \tanh\sigma \sinh 2\psi \right) (\sinh^2 \psi + \tanh^2 \sigma)^{-1}$$
$$- \lambda \sinh\psi \cosh\psi (\sinh^2 \psi + \tanh^2 \sigma)^{-1/2}, \quad (12)$$

$$y_3 = \left( -\psi^* \tanh\sigma \cosh 2\psi + \sigma^* \operatorname{sech}^2 \sigma \sinh\psi \cosh\psi \right) (\sinh^2 \psi + \tanh^2 \sigma)^{-1}$$
$$+ \lambda \sinh^2 \psi (\sinh^2 \psi + \tanh^2 \sigma)^{-1/2}. \quad (13)$$

If we put $\sigma = 0$, $\sigma^* \neq 0$ in equations (11), (12) and (13), then we get

$$y_1 = \psi^*, \quad (14)$$

$$y_2 = -\sigma^* - \lambda \cosh\psi, \quad (15)$$

$$y_3 = \sigma^* \frac{\cosh\psi}{\sinh\psi} + \lambda \sinh\psi. \quad (16)$$

Equations (14), (15) and (16) give a two-parameter family (linear congruence) of the timelike straight lines which are the intersection of the planes $y_1 = \psi^*$ and the timelike

ruled surfaces given by

$$\cosh^2\psi\left(y_2+\frac{\lambda}{\cosh\psi}\right)^2-y_3^2\sinh^2\psi=0. \qquad (17)$$

Thus we give the following theorem:

***Theorem 3.1.*** *During the dual hyperbolic conchoidal motion, in the case of $\sigma=0$, $\sigma^*\neq 0$, the Study map in $\mathrm{IR}_1^3$ of the orbit which is drawn on the $H'$ by $\tilde{X}=\tilde{v}_1$ are the straight lines which are the intersections of the planes $y_1=\psi^*$ and the timelike ruled surfaces given by*

$$\cosh^2\psi\left(y_2+\frac{\lambda}{\cosh\psi}\right)^2-y_3^2\sinh^2\psi=0.$$

Now, let us take $P=p+\varepsilon p^*=0$ in equation (4). In this case, $\tilde{X}=\tilde{v}_3$, thus, from equations (8), (9) and (10)

$$y_1=\psi^*, \qquad (18)$$

$$y_2=\lambda\sinh\psi, \qquad (19)$$

$$y_3=\lambda\cosh\psi, \qquad (20)$$

are obtained. From equations (18), (19) and (20) we have

$$\left.\begin{array}{l}y_3^2-y_2^2=\lambda^2,\\ y_1=\psi^*.\end{array}\right\} \qquad (21)$$

Thus we have the following theorem.

***Theorem 3.2.*** *During the dual hyperbolic conchoidal motion $H/H'$, in the case of $P=p+\varepsilon p^*=0$ in equation (4), the Study map of the orbit which is drawn on the $H'$ by $\tilde{X}=\tilde{v}_3$ is the congruence,*

$$y_3^2-y_2^2=\lambda^2,$$
$$y_1=\psi^*.$$

Let us now give the analysis of the orbit of $\tilde{v}_2$ during the dual hyperbolic conchoidal motion. We know that

$$\tilde{v}_2 = \left[ -\tilde{\psi}\left(\sinh^2\tilde{\psi} + \tanh^2\tilde{\varDelta}\right)^{-1/2}, \ -\cosh\tilde{\psi}\ \tanh\tilde{\varDelta}\left(\sinh^2\tilde{\psi} + \tanh^2\tilde{\varDelta}\right)^{-1/2}, \right.$$
$$\left. -\sinh\tilde{\psi}\ \tanh\tilde{\varDelta}\left(\sinh^2\tilde{\psi} + \tanh^2\tilde{\varDelta}\right)^{-1/2} \right]. \tag{22}$$

From equation (22), we obtain

$$v_2 = \left(-A^{-1/2}\sinh\psi, \ -A^{-1/2}\tanh\sigma\cosh\psi, \ -A^{-1/2}\tanh\sigma\sinh\psi\right) \tag{23}$$

and

$$\begin{aligned}v_2^* = (&\psi^* A^{-3/2}\sinh^2\psi\cosh\psi - \psi^* A^{-1/2}\cosh\psi + \sigma^* A^{-3/2}\tanh\sigma\,\text{sech}^2\sigma\sinh\psi,\\ &-\sigma^* A^{-1/2}\text{sech}^2\sigma\cosh\psi + \psi^* A^{-3/2}\tanh\sigma\sinh\psi\cosh^2\psi\\ &-\psi^* A^{-1/2}\tanh\sigma\sinh\psi + \sigma^* A^{-3/2}\tanh^2\sigma\,\text{sech}^2\sigma\cosh\psi,\\ &+\psi^* A^{-3/2}\tanh\sigma\sinh^2\psi\cosh\psi - \psi^* A^{-1/2}\tanh\sigma\cosh\psi\\ &+\sigma^* A^{-3/2}\tanh^2\sigma\,\text{sech}^2\sigma\sinh\psi - \sigma^* A^{-1/2}\text{sech}^2\sigma\sinh\psi)\end{aligned} \tag{24}$$

where $v_2$ and $v_2^*$ are the real and dual parts of $\tilde{v}_2$, respectively, and $A = \sinh^2\psi + \tanh^2\sigma$.

Equations (23) and (24) depend on two parameters $\psi$ and $\psi^*$ so equations (23) and (24) represent a timelike line congruence in $IR_1^3$.

Let $g$ denote the position vector of an arbitrary point $G(g_1, g_2, g_3)$ of a fixed timelike line $x$ of the timelike congruence in $IR_1^3$. Then considering equation (7) we have

$$g = v_2(\psi, \psi^*) \times v_2^*(\psi, \psi^*) + u v_2(\psi, \psi^*). \tag{25}$$

Since $(g_1, g_2, g_3)$ are the coordinates of $G$ we have

$$g_1 = A^{-1}\psi^* \tanh^2\sigma + uA^{-1/2}\sinh\psi, \tag{26}$$

$$g_2 = -A^{-1}\sigma^* \tanh^2\sigma + uA^{-1/2}\tanh\sigma\cosh\psi, \tag{27}$$

$$g_3 = -A^{-1}\sigma^* \text{sech}^2\sigma\sinh\psi\cosh\psi + A^{-1}\psi^* \tanh\sigma + uA^{-1/2}\tanh\sigma\sinh\psi, \tag{28}$$

where $A = \sinh^2\psi + \tanh^2\sigma$.

If $\sigma = 0$, $\sigma^* \neq 0$ in equations (26), (27) and (28), then

$$g_1 = u, \ g_2 = -\sigma^*, \ g_3 = -\sigma^*\coth\psi \tag{29}$$

are obtained. In this case, if we choose $u = k\psi$ ($k$ constant) in equation (29) we have

$$g_1 = k \tanh^{-1} \frac{g_2}{g_3}, \tag{30}$$

which is a Lorentzian helicoid.

If $\varphi = 0$, $\varphi^* \neq 0$ in equations (26), (27) and (28), then

$$g_1 = \varphi^*, \quad g_2 = u, \quad g_3 = \varphi^* \tanh \sigma \tag{31}$$

are obtained. If we choose $u = k\sigma$ ($k$ constant) in equation (31) we have

$$g_2 = k \tanh^{-1} \frac{g_3}{g_1},$$

which is also a Lorentzian helicoid.